\documentclass{amsart}
\usepackage{amsmath}
\usepackage{amsthm}
\usepackage{amssymb}
\usepackage{epsfig}

\newcommand{\hs}{{\mathcal{H}(m,n)}}
\newcommand{\hso}{{\mathcal{H}^o(m,n)}}
\newcommand{\halfhs}{{H(m,n)}}

\newcommand{\N}{{N(kn,(k+1)n)}}

\theoremstyle{plain}
\newtheorem{thm}{Theorem}
\newtheorem{cor}{Corollary}

\title{Tiling Parity Results and the Holey Square Solution}
\author{Bridget Eileen Tenner}
\date{October 7, 2004}

\begin{document}

\begin{abstract}
We prove combinatorially that the parity of the number of domino tilings of a region is equal to the parity of the number of domino tilings of a particular subregion. Using this result we can resolve the holey square conjecture.  We furthermore give combinatorial proofs of several other tiling parity results, including that the number of domino tilings of a particular family of rectangles is always odd.
\end{abstract}

\maketitle

\section{Introduction}

The number of domino tilings of the $2n \times 2n$ square with a centered hole of size $2m \times 2m$, a figure known as the \emph{holey square} and denoted $\hs$, was conjectured by Edward Early to have the form $2^{n-m} (2k_{m,n}+1)^2$.  Although the conjecture has remained unsolved until now, specific cases were known, for example see \cite{tauraso}.  In this paper, we answer the general conjecture in the affirmative, primarily via a theorem about tiling parity that has applications beyond the problem of the holey square.  We also give combinatorial meaning to the odd factor $2k_{m,n}+1$ in Early's conjecture, and demonstrate other consequences of the parity theorem.

\begin{figure}[htbp]
\centering
\epsfig{file=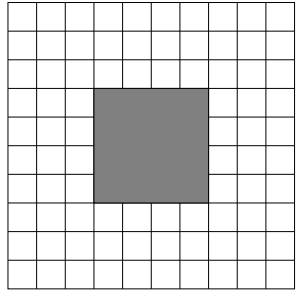}
\caption{The holey square $\mathcal{H}(2,5)$.  Throughout this paper, shading indicates a portion of the figure that is excluded from the region.}
\end{figure}

As this paper solely concerns domino tilings, all tilings discussed can be assumed to be domino tilings.  Following Pachter's notation in \cite{pachter}, we write $\# R$ for the number of tilings of the region $R$, and $\#_2 R$ for the parity of the number of tilings of $R$.  When we are only concerned with the configuration of part of a region, we may only draw this portion, indicating that the undrawn portion is arbitrary.

\section{A tiling parity result}

\begin{figure}[htbp]\label{stp}
\centering
\epsfig{file=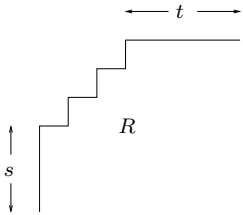}
\caption{An $(\{s,t\};4)$-corner.}
\end{figure}

In this section, we present a theorem regarding the parity of the number of domino tilings of a region.  This result depends only on a very local property of the region, and makes no further assumptions regarding symmetry, planarity, or any other aspect of the region.  Say that a region $R$ has an \emph{$(\{s,t\}; 1)$-corner} if there is a convex corner in $R$ where the segments bounding this corner have lengths $s$ and $t$.  For $p > 1$ and $\min\{s,t\} \ge 2$, say that $R$ has an \emph{$(\{s,t\}; p)$-corner} if there is a $(\{1,s\};1)$-corner, a $(\{1,t\};1)$-corner, and $p-2$ distinct $(\{1,1\};1)$-corners configured as in Figure~\ref{stp}.

If one of the segments, say the segment of length $s$, in an $(\{s,t\};p)$-corner forms an $(\{s,t'\}; p')$-corner at its other endpoint for some $t'$ and $p'$, say that each of these corners is \emph{walled} at $s$.

Let an \emph{$(\{i,j\};p)$-strip} be a subregion of $i+j+2p-3$ squares that has an $(\{i,j\};p)$-corner.  For $2 \le k \le \min\{s,t\}$, say that an $(\{s,t\};p)$-corner in a region $R$ is \emph{complete up to $k$} if that corner is $i$-complete for $i = 2, \ldots, k$, where an $(\{s,t\};p)$-corner is \emph{$2$-complete} if $2 \le \min\{s,t\}$, and \emph{$i$-complete} for $2 < i \le \min\{s,t\}$ if the following inductive definition is true.

\begin{enumerate}
\item Let $C$ be the $(\{i,i\}; p)$-strip in the $(\{s,t\};p)$-corner.  Let $x$ and $y$ be the two squares adjacent to the ends of $C$ but not along the edges forming the $(\{s,t\};p)$-corner.  If either $x$ or $y$ is in $R$, then the $(\{i-1,i-1\};p)$-strip between $x$ and $y$, inclusively, all of whose squares are adjacent to $C$, must also be a subregion of $R$.

\begin{figure}[htbp]
\centering
\epsfig{file=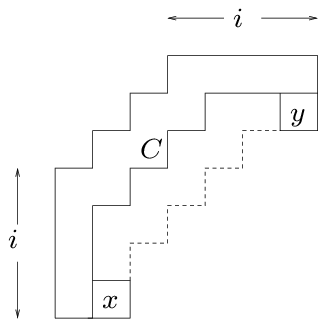}
\end{figure}

\item Consider the $(\{s',t'\};p)$-corner formed by removing $C$ from $R$.  If $2 \le i-2 \le \min\{s',t'\}$, then this corner must be complete up to $i-2$.

\item The subregion of squares examined at each step of this definition must be the dual of a grid graph.
\end{enumerate}

For example, if $R$ has an $(\{s,t\};1)$-corner that is complete up to $3$, then this corner must have one of the following forms.
\begin{equation*}
\parbox{.6in}{\epsfig{file=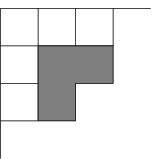}}
\hspace{.2in}
\parbox{.6in}{\epsfig{file=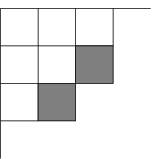}}
\hspace{.2in}
\parbox{.6in}{\epsfig{file=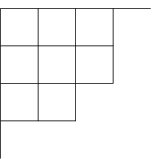}}
\end{equation*}

\noindent Similarly, if $R$ has an $(\{s,t\};2)$-corner that is complete up to $3$, the possibilities for this corner are as follows.
\begin{equation*}
\parbox{.9in}{\epsfig{file=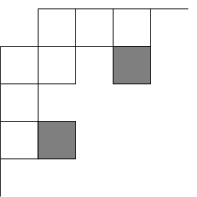}}
\hspace{.2in}
\parbox{.9in}{\epsfig{file=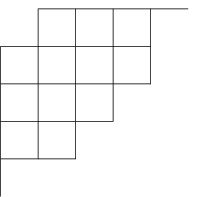}}
\end{equation*}

The different tilings of $R$ can be categorized by the manner in which an $(\{s,t\};p)$-corner is tiled.  For example, if $R$ has an $(\{s,t\};4)$-corner, then
\begin{equation*}
\# R = \# \hspace{.05in}
\parbox{1.2in}{\epsfig{file=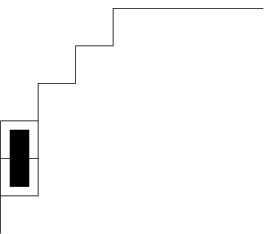}}
+ \# \hspace{.05in}
\parbox{1.2in}{\epsfig{file=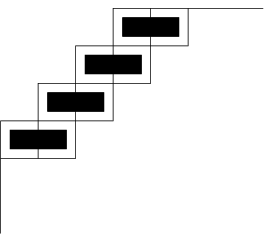}}.
\end{equation*}

\noindent Suppose $R$ has an $(\{s,t\};1)$-corner where $2 \le \min\{s,t\}$.  Then
\begin{equation*}
\# R = \# \hspace{.05in}
\parbox{.6in}{\epsfig{file=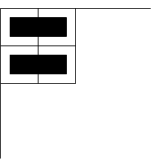}}
+ \# \hspace{.05in}
\parbox{.6in}{\epsfig{file=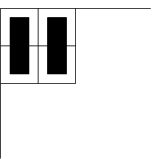}}
+ \# \hspace{.05in}
\parbox{.6in}{\epsfig{file=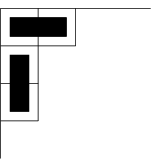}}
+ \# \hspace{.05in}
\parbox{.6in}{\epsfig{file=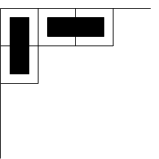}},
\end{equation*}

\noindent where the corner drawn in each of the above figures is the particular $(\{s,t\};1)$-corner in $R$.  If $R$ does not include the entire region tiled in one of these figures, that term is zero.  Notice that the first two of these figures tile the same subregion of $R$, so in fact,
\begin{equation*}
\#_2 R = \#_2 \hspace{.05in}
\parbox{.6in}{\epsfig{file=cornerhv.eps}}
+ \#_2 \hspace{.05in}
\parbox{.6in}{\epsfig{file=cornervh.eps}}.
\end{equation*}

\begin{thm}[Parity Theorem]

Suppose that a region $R$ has an $(\{s,t\};p)$-corner.  For any $2 \le k \le \min\{s,t\}$, if this corner is complete up to $k$, then
\begin{equation}\label{open}
\#_2 R = \#_2 \hspace{.05in}
\parbox{1.5in}{\epsfig{file=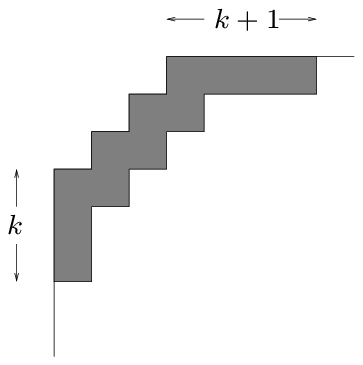}}
+ \#_2  \hspace{.05in}
\parbox{1.6in}{\epsfig{file=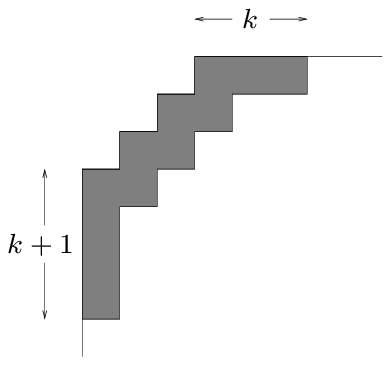}}.
\end{equation}

\noindent If $p = 1$, then \eqref{open} also holds for $k=1$.  Furthermore, for any $p$, if $s \le t$, the corner is complete up to $s$, and the corner is walled at $s$, then 
\begin{equation}\label{wall} 
\#_2 R = \#_2  \hspace{.05in}
\parbox{1.5in}{\epsfig{file=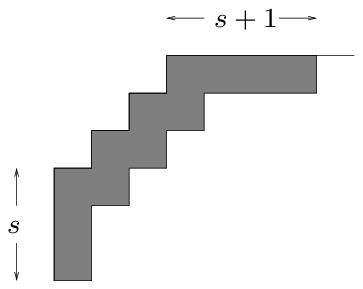}}.
\end{equation}

\end{thm}

\begin{proof}

We prove the theorem by first inductively showing that for a particular $p$, the statement is true for all $(\{s,t\};p)$-corners complete up to $k \le \min\{s,t\}$.  Then we induct on $p$.  For any $p$, if \eqref{open} holds for a particular $k$, then certainly \eqref{wall} holds if there is a wall at $s=k$, since one of the tilings pictured in the statement of \eqref{open} is impossible because of the wall at $s$, so this term is zero.

Suppose first that $p=1$.  That \eqref{open} holds for $k=1$ is trivial, and the case $k=2$ was shown above.  Suppose that the theorem holds for all $k < K \le \min\{s,t\}$.  Suppose that an $(\{s,t\};1)$-corner of $R$ is complete up to $K$.  This corner must also be complete up to $K-1$, so we can apply \eqref{open} for $k=K-1$.
\begin{equation}\label{K-1} 
\#_2 R = \#_2  \hspace{.05in}
\parbox{1.1in}{\epsfig{file=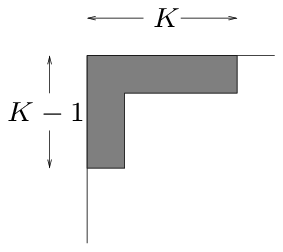}}
+ \#_2 \hspace{.05in}
\parbox{1in}{\epsfig{file=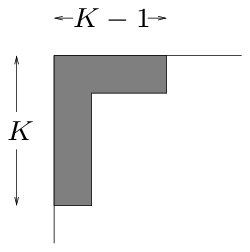}}
\end{equation}

Consider the first of these regions, and look at the possible ways to tile the square next to the shorter leg of the removed region.
\begin{equation*} 
\#_2  \hspace{.05in}
\parbox{1.1in}{\epsfig{file=K-1K.eps}}
= \#_2 \hspace{.05in}
\parbox{1.1in}{\epsfig{file=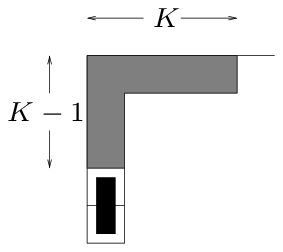}}
+ \#_2 \hspace{.05in}
\parbox{1.1in}{\epsfig{file=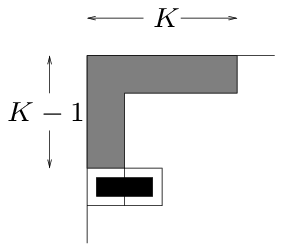}}
\end{equation*}

\noindent The latter of these creates a $(\{K-2,K-1\};1)$-corner walled at $K-2$, and this corner is complete up to $K-2$ because of the inductive definition of $K$-completeness.  Thus we can apply \eqref{wall} for $k = K-2$ to get 
\begin{equation}\label{doubleterm} 
\#_2  \hspace{.05in}
\parbox{1.1in}{\epsfig{file=K-1K.eps}}
= \#_2 \hspace{.05in}
\parbox{1.1in}{\epsfig{file=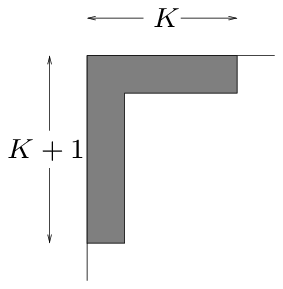}}
+ \#_2 \hspace{.05in}
\parbox{1.2in}{\epsfig{file=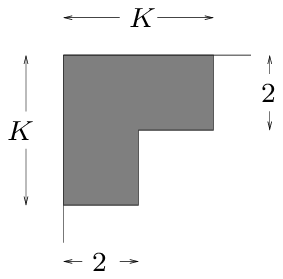}}.
\end{equation}

There is an analogous equation for the second region on the right-hand side of \eqref{K-1}, and the last term in \eqref{doubleterm} also appears in this.  Since we are considering parity, these terms cancel, leaving
\begin{equation*}
\#_2 R = \#_2  \hspace{.05in}
\parbox{1.1in}{\epsfig{file=goodK+1K.eps}}
+ \#_2 \hspace{.05in}
\parbox{1.2in}{\epsfig{file=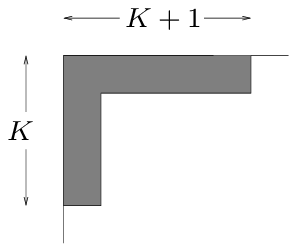}}
\end{equation*}

\noindent which proves \eqref{open} for $k = K$, and consequently \eqref{wall} for $k = K$.  Thus the result holds for all $(\{s,t\};1)$-corners.

Now suppose that the theorem is true for all $1 \le p < P$.  Suppose that a region $R$ has a $(\{s,t\};P)$-corner where $2 \le \min\{s,t\}$.  First we must show that the result holds for this corner for $k=2$, and then we can induct on $k$.

Consider the $(\{s,1\};1)$-corner in this $(\{s,t\};P)$-corner.  We can either place a vertical domino or a horizontal domino in this corner, which gives the following.
\begin{equation*}
\# R = \# \hspace{.05in}
\parbox{1.1in}{\epsfig{file=pv.eps}}
+ \# \hspace{.05in}
\parbox{1.1in}{\epsfig{file=ph.eps}}
\end{equation*}

\noindent Observe that the first of these possibilities creates a $(\{3,t\}; P-1)$-corner which is necessarily complete up to $2$.  The second possible tiling creates a $(\{2,s-1\};1)$-corner.  To the former, apply the theorem for $p=P-1$ and $k=2$, and apply the theorem for $p=1$ and $k=1$ to the latter.  Both of these results are already known by the induction hypothesis.  Thus 
\begin{eqnarray*}
\# R &=&\# \hspace{.05in}
\parbox{1.1in}{\epsfig{file=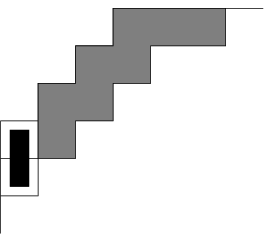}}
+ \# \hspace{.05in}
\parbox{1.1in}{\epsfig{file=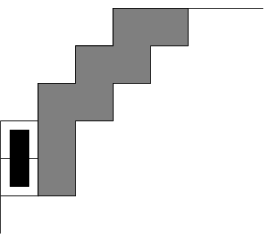}}\\
& &+ \# \hspace{.05in}
\parbox{1.1in}{\epsfig{file=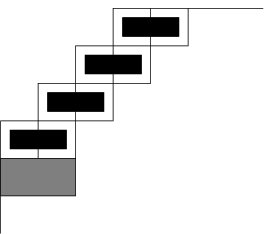}}
+ \# \hspace{.05in}
\parbox{1.1in}{\epsfig{file=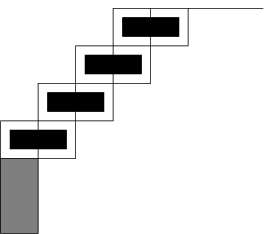}}.
\end{eqnarray*}

\noindent The second and third of these tile the same region, so when we consider $\#_2 R$, we see that the theorem is true for $k=2$.

Now suppose that for this $P$, the theorem holds for all $2 \le k < K \le \min\{s,t\}$, and that there is an $(\{s,t\};P)$-corner of $R$ complete up to $K$.  As in the case of a $(\{s,t\};1)$-corner, this corner must also be complete up to $K-1$, so we can apply \eqref{open} for $k=K-1$.  The remainder of the proof follows from methods that are entirely analogous to the case when $p=1$.  This completes the proof for all $(\{s,t\};p)$-corners complete up to $k \in [2, \min\{s,t\}]$.

\end{proof}

\section{The holey square}

Counting the number of tilings of the holey square can be simplified using two known results, both of which rely on the holey square's symmetries.  A result of Jockusch, see \cite{jockusch}, indicates that $\# \hs$ is either a square or twice a square.  In \cite{ciucu}, Ciucu defines a Klein-symmetric graph as a bipartite graph with $180$-degree rotational symmetry and a reflective axis of symmetry that separates the graph.  Say that a region in the plane is Klein-symmetric if its dual is a Klein-symmetric graph.  The region $\hs$ is Klein-symmetric with each diagonal serving as a reflexive axis of symmetry, so we can apply the factorization theorem in \cite{ciucu} to conclude that $\# \hs = 2^{n-m} (\# \halfhs)^2$ for a region $\halfhs$ defined as follows.  Consider the $2n \times 2n$ square, coordinatized so that the lower left corner is at the origin and the upper right corner has coordinates $(2n,2n)$.  Divide the square into two congruent halves by the two-unit segments

\begin{center}
$\big\{[2t,2t+2] \times \{2n - (2t+1)\} : t = 0, \ldots, n-1\big\}$

$\cup \big\{\{2n-2t\} \times [2t-1,2t+1] : t = 1, \ldots, n-1\big\}.$
\end{center}

\begin{figure}[htbp]\label{halfhs}
\centering
\epsfig{file=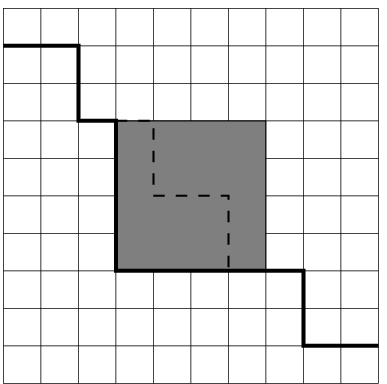}
\caption{The region below the heavy line is $H(2,5)$.}
\end{figure}

\noindent Now remove the center $2m \times 2m$ square from the region.  This leaves two congruent regions, denote each by $\halfhs$.  Notice that in any tiling of $\halfhs$, one particular domino must always be in place.  In Figure~\ref{halfhs}, this is the domino occupying the two bottom-rightmost squares.  Let $H'(m,n)$ be $\halfhs$ with these two squares removed.  There is an obvious bijection between tilings of $\halfhs$ and tilings of $H'(m,n)$, so $\# \halfhs = \# H'(m,n)$, and in particular these numbers have the same parity.

To prove the conjecture we need to prove that $\#_2 \halfhs = 1$.   Our proof is inductive for fixed $m$, with base case $n = m+1$.  This base case is trivial, as it is easy to see that $\#_2 H(m,m+1) = 1$ and $\mathcal{H}(m,m+1)$ can be tiled in two ways.  

\begin{cor}

For all $m$ and $n > m$, $\# \hs = 2^{n-m}(2k_{m,n}+1)^2$, where the factor $2k_{m,n}+1$ is equal to $\# \halfhs$.

\end{cor}

\begin{proof}

\begin{sloppypar}

Assume that $\#_2 \halfhs = 1$.  As discussed above, it remains to show only that $\#_2 H(m,n+1) = 1$.  Consider the region $H'(m,n+1)$ which has a \mbox{$(\{2n,2n+1\};1)$}-corner walled at $2n$, and this corner is complete up to $2n$.  Therefore we can apply the parity theorem to this corner, specifically \eqref{wall}.  The subregion indicated by the right-hand side of \eqref{wall} is actually $\halfhs$ reflected across the line $y=x$.  Therefore $\#_2 H'(m,n+1) = \#_2 \halfhs = 1$.  This completes the proof since $\#_2 H'(m,n+1) = \#_2 H(m,n+1)$, answering affirmatively the question posed by Early, and giving a combinatorial meaning to the odd factor in $\# \hs$.

\end{sloppypar}

\end{proof}

Analogous to $\hs$, let $\hso$ be the $(2n+1) \times (2n+1)$ square with a centered hole of size $(2m+1)\times (2m+1)$.

\begin{cor}

For all $m$ and $n > m$, $\# \hso = 2^{n-m}(2k_{m,n}'+1)^2$.

\end{cor}

\begin{proof}

The proof is analogous to the proof of the previous corollary, and once again the odd factor $2k_{m,n}'+1$ is the number of domino tilings of a particular region.

\end{proof}

\section{Further applications of the parity theorem}

In addition to determining the number of domino tilings of the holey square, the parity theorem can be applied to other regions.  One easy consequence is the following.

\begin{cor}\label{samelength}

Suppose a region $R$ has an $(\{s,s\};p)$-corner that is complete up to $s$ and walled at $s$ along both sides.  Then $\# R$ is even.

\end{cor}

\begin{proof}

Much like in the proof of \eqref{wall}, notice that both of the tilings depicted in \eqref{open} are impossible, so both terms on the right side of the equation are zero.

\end{proof}

For another application, let $N(a,b)$ be the $a \times b$ rectangle, and notice that each corner is complete up to $\min\{a,b\}$.  We can repeatedly apply the parity theorem to $\N$ to obtain the following result.

\begin{cor}\label{kncor}

For all positive integers $k$ and $n$, $\# \N$ is odd.

\end{cor}

\begin{proof}

\begin{sloppypar}

To describe the process more precisely, suppose that $\N$ is oriented so that the sides of length $kn$ are vertical.  The upper left corner is complete up to $kn$, so apply the parity theorem to remove the \mbox{$(\{kn,kn+1\};1)$}-strip from this corner.  Similarly, if the values of $k$ and $n$ are sufficiently large, we can remove the \mbox{$(\{n-1,n\};1)$}-strip from the upper right corner, the \mbox{$(\{(k-1)n,(k-1)n+1\};1)$}-strip from the lower right corner, and the \mbox{$(\{2n-2,2n-1\};1)$}-strip from the lower left corner (where each of these corners are of subsequent subregions of $\N$).

\end{sloppypar}

\begin{figure}[htbp]
\centering
\epsfig{file=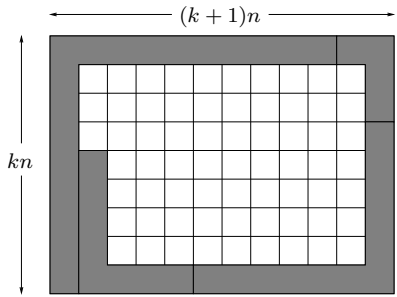}
\caption{$\N$ with the first four removed strips marked, as implied by the parity theorem.}
\end{figure}

The process of removing each strip can be summarized in the following table describing how many squares are removed from the sides of the region, starting with $\N$, after each application of the parity theorem.

\begin{center}
{\small
\begin{tabular}{c|l|l|l|l}
App. & Left & Top & Right & Bottom\\
\hline
$1$ & $kn$ & $kn + 1$ & --- & ---\\
$2$ & --- & $n-1$ & $n$ & --- \\
$3$ & --- & --- & $(k-1)n$ & $(k-1)n+1$ \\
$4$ & $2n-1$ & --- & --- & $2n-2$ \\
$\vdots$ & $\ \vdots$ & $\ \vdots$ & \ $\vdots$ & $\ \vdots$ \\
$4j-3$ & $(k-2j+2)n$ & $(k-2j+2)n+1$ & --- & --- \\
$4j-2$ & --- & $(2j-1)(n-1)$ & $(2j-1)(n-1) + 1$ & --- \\
$4j-1$ & --- & --- & $(k - 2j+1)n$ & $(k - 2j+1)n + 1$ \\
$4j$ & $2j(n-1) + 1$ & --- & --- & $2j(n-1)$
\end{tabular}
}
\end{center}

We can continue to apply the parity theorem until the $(\{i,i+1\};1)$-strip we remove is actually from an $(\{i,i+1\};1)$-corner.  To determine when this might happen, we need to solve any of the following equations, where the term subtracted from the left side of each refers to the squares occupied by previously removed corners.
\begin{equation}\label{finishleft}
kn - (2j - 1) = 2j(n - 1) + 1
\end{equation}
\begin{equation}\label{finishtop}
(k+1)n - (2j - 2) = (k-2j + 2)n + 1
\end{equation}
\begin{equation}\label{finishright}
kn - (2j - 2) = (2j-1)(n-1) + 1
\end{equation}
\begin{equation}\label{finishbottom}
(k+1)n - (2j - 1) = (k-2j+1)n + 1
\end{equation}

These equations correspond to the final $(\{i,i+1\};1)$-strip being removed from the lower left corner, the upper left corner, the upper right corner, and the lower right corner, respectively.  The solution to \eqref{finishleft} is $n=0$ or $k = 2j$, the solution to \eqref{finishtop} is $n = 1$ or $2j-1 = 0$, the solution to \eqref{finishright} is $n=0$ or $k = 2j-1$, and the solution to \eqref{finishbottom} is $n = 1$ or $2j = 0$.  Since $n$, $k$, and $j$ are positive integers, several of these cases are impossible.  Thus the applications of the parity theorem cease in the manner described above when $n=1$ or when $2j = k$ or $2j-1 = k$, depending on the parity of $k$.

Suppose that $n > 1$, and let $k$ be even.  Set $j$ to be $k/2$.  After removing the $(2k)^{\text{th}}$ strip from $\N$, it is not hard to see that we are left with the subregion of $\N$ formed by removing the top $j$ rows, the bottom $j$ rows, the left $j+1$ columns, and the right $j$ columns.  The remaining region is $N(k(n-1), (k+1)(n-1))$, so 
\begin{equation}\label{rectreduction}
\#_2 \N = \#_2 N(k(n-1), (k+1)(n-1)).
\end{equation}

\noindent If, on the other hand, $k$ is odd, let $j = (k+1)/2$.  After removing the $(2k)^{\text{th}}$ strip from $\N$, we likewise find that the resulting subregion of $\N$ is $\N$ with the top $j$ rows, the bottom $j-1$ rows, the left $j$ columns, and the right $j$ columns removed.  Since $k = 2j-1$, this once again gives \eqref{rectreduction}.

Therefore, $\#_2 \N = \#_2 N(k,k+1)$ for all positive integers $k$ and $n$.  Applying the parity theorem once to any corner of $N(k,k+1)$ indicates that
\begin{equation*}
\#_2 N(k,k+1) = \#_2 N(k-1,k)
\end{equation*}

\noindent for all $k > 1$.  Therefore $\#_2 \N = \#_2 N(1,2) = 1$ for all positive integers $k$ and $n$.

\end{proof}

It should be noted that there are other ways to obtain this result, for example using Kasteleyn's formula or determinants of particular matrices.  However, these methods tend to be much more analytic, and thus somewhat less intuitively clear, than the combinatorial proof presented here.

\begin{sloppypar}

In Corollary~\ref{kncor}, we studied $(\{s,t\};1)$-corners.  We conclude this section by considering more general types of corners.  Let $T(i,j,p)$ be the region with \mbox{$i+p-1$} rows, whose rows from top to bottom consist of the following number of squares: $j$, \mbox{$j+2$}, $\ldots$, \mbox{$j+2(p-1)$}, $\ldots$, \mbox{$j+2(p-1)$}, where there are $i$ rows of \mbox{$j+2(p-1)$} boxes.  Similarly, let $D(i,j,p)$ be the region with $i+2(p-1)$ rows, whose rows from top to bottom consist of the following number of squares: $j$, \mbox{$j+2$}, $\ldots$, \mbox{$j+2(p-1)$}, $\ldots$, \mbox{$j+2(p-1)$}, \mbox{$j+2(p-1)-2$}, $\ldots$, \mbox{$j+2$}, $j$, where there are $i$ rows of \mbox{$j+2(p-1)$} boxes.  In each of these regions, the centers of the rows are aligned along a vertical line.

\end{sloppypar}

\begin{figure}[htbp]
\centering
$\begin{array}{c@{\hspace{.5in}}c}
\multicolumn{1}{l}{\mbox{\bf{(a)}}} & \multicolumn{1}{l}{\mbox{\bf{(b)}}}\\
[-.25cm]
\epsfig{file=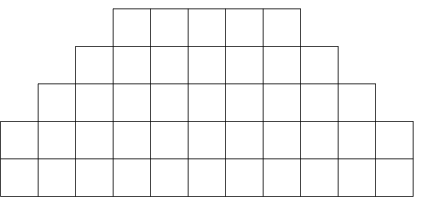} &\epsfig{file=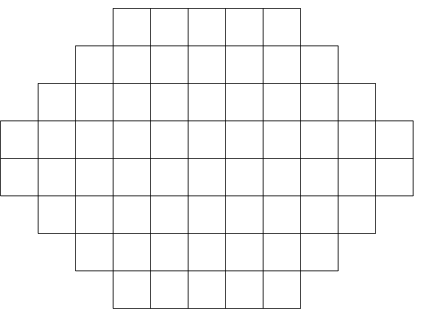}
\end{array}$
\caption{$\bf{(a)}$ The region $T(2,5,4)$.  $\bf{(b)}$ The region $D(2,5,4)$.}
\end{figure}

We can use the parity theorem to determine the parities of many different $T(i,j,p)$ and $D(i,j,p)$.  Examples of such applications are given below.  Also notice that for $T(i,j,p)$ to have an even number of squares (and hence be possibly tilable), either $j$ or $i+p-1$ must be even.  Similarly, for $D(i,j,p)$ to have an even number of squares, either $j$ or $i$ must be even.

\begin{cor}

\begin{enumerate}

\setlength{\itemsep}{2mm}

\item $\#_2 T(k,k,p) = 0$.

\item $\#_2 T(k,k+1,p) = \left\{ \begin{array} {c@{\quad:\quad}l}
0 & p > 1 \text{ and } k \text{ is even}; \\ 1 & p=1 \text{ or } k \text{ is odd}. \end{array} \right.$

\item $\#_2 T(k,k+2,p) = \left\{ \begin{array} {c@{\quad:\quad}l}
0 & k \text{ is odd}; \\ 1 & k \text{ is even}. \end{array} \right.$

\item $\#_2 T(k,2k-1,p) = 0$.

\item $\#_2 T(k,2k,p) = 1$.

\item $\#_2 T(k,2k+1,p) = 0$.

\item $\#_2 T(k,2k+2,p) = 1$.

\end{enumerate}
\end{cor}

\begin{proof}

Each of these follows inductively from (sometimes several) straightforward applications of the parity theorem. Outlines of the proofs are as follows.

\begin{enumerate}

\item This follows immediately from Corollary~\ref{samelength}.
\begin{equation*}
\epsfig{file=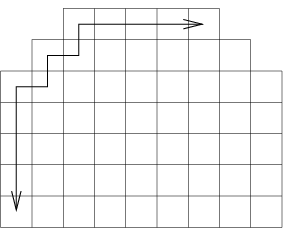}
\end{equation*}

\item First notice that $T(k,k+1,1)$ is a $k \times (k+1)$ rectangle, so the result for $p=1$ follows from Corollary~\ref{kncor}.  For $p>1$, it is not hard to see that $\#_2 T(k,k+1,p) = \#_2 T(k-1,k+1,p-1)$, so for $k>1$ we must use the answer from the subsequent part.  Also, $T(1,2,p)$ has one tiling.
  \begin{equation*}
\epsfig{file=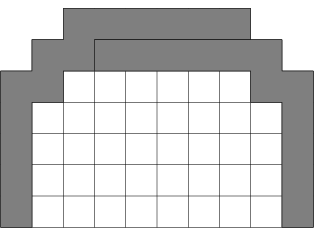}
\end{equation*}

\item We see that $\#_2 T(k,k+2,p) = \#_2 T(k-2,k,p)$, and $\#_2 T(1,3,p) = 0$ for all $p$, while $\#_2 T(2,4,p) = \#_2 T(1,2,p) = 1$ for all $p$.
\begin{equation*}
\epsfig{file=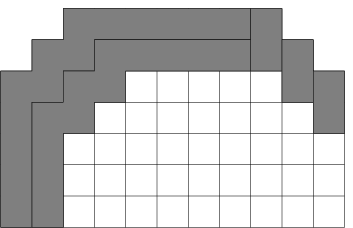}
\end{equation*}

\item Several applications of the parity theorem imply that $\#_2 T(k,2k-1,p) = \#_2 T(k-2, 2(k-2)-1, p)$, and for both $k = 1$ and $k=2$ this value is $0$.
\begin{equation*}
\epsfig{file=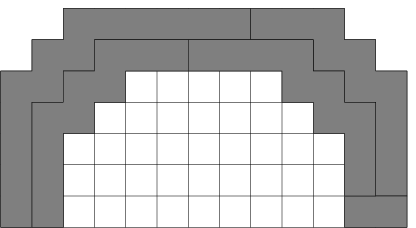}
\end{equation*}

\item $\#_2 T(k,2k,p) = \#_2 T(k-1, 2(k-1), p)$ and $\#_2 T(1,2,p) = 1$.
\begin{equation*}
\epsfig{file=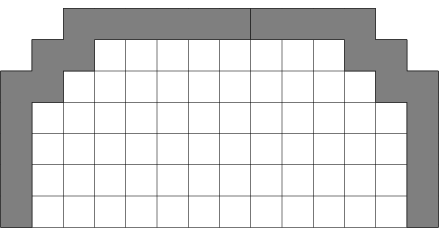}
\end{equation*}

\item After applying the parity theorem once, we can apply Corollary~\ref{samelength}.
\begin{equation*}
\epsfig{file=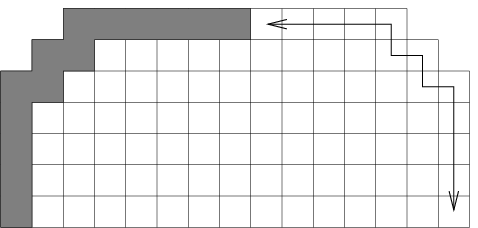}
\end{equation*}

\item Two applications of the parity theorem imply that $\#_2 T(k,2k+2,p) = \#_2 T(k-1, 2k, p)$, and $\#_2 T(1,4,p) = 1$ for all $p$.
\begin{equation*}
\epsfig{file=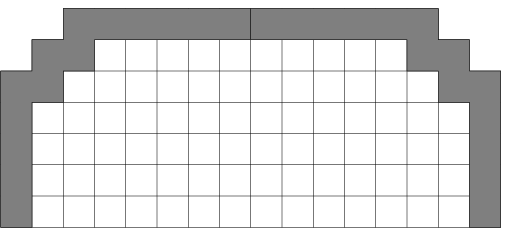}
\end{equation*}

\end{enumerate}

\end{proof}

\begin{cor}

\begin{enumerate}

\setlength{\itemsep}{2mm}

\item $\#_2 D(k,k,p) = 0$.

\item $\#_2 D(k,k+1,p) = 1$.

\item $\#_2 D(k,k+2,p) = 0$.

\item $\#_2 D(k,2k-1,p) = \left\{ \begin{array} {c@{\quad:\quad}l}
0 & k \text{ is odd}; \\ 1 & k \text{ is even}. \end{array} \right.$

\item $\#_2 D(k,2k+1,p) = 0$.


\end{enumerate}
\end{cor}

\begin{proof}

As before, these results follow inductively from the parity theorem.

\begin{enumerate}

\item This follows immediately from Corollary~\ref{samelength}.
\begin{equation*}
\epsfig{file=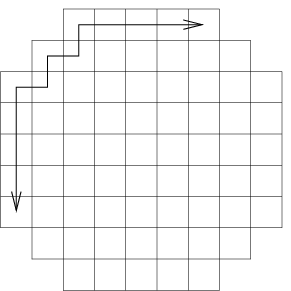}
\end{equation*}

\item After applying the parity theorem twice, we see that $\#_2 D(k,k+1,p) = \#_2 D(k,k+1,p-1)$, and $\#_2 D(k,k+1,1) = 1$.
\begin{equation*}
\epsfig{file=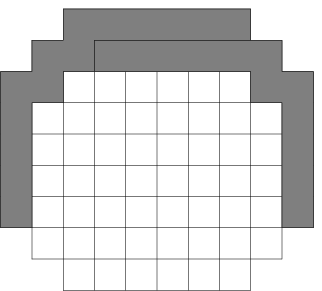}
\end{equation*}

\item If $p > 1$, then this follows from Corollary~\ref{samelength} after one application of the parity theorem.  If $p=1$, then $\#_2 D(k,k+2,1) = \#_2 D(k-2,k,1)$, and this is $0$ for $k=1$ and $k=2$.
\begin{equation*}
\epsfig{file=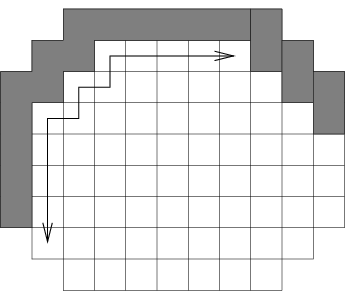}
\end{equation*}

\item We see that $\#_2 D(k,2k-1,p) = \#_2 D(k-2, 2(k-2)-1,p)$, and $\#_2 D(1,1,p) = 0$ for all $p$, while $\#_2 D(2,3,p) = 1$ for all $p$.
\begin{equation*}
\epsfig{file=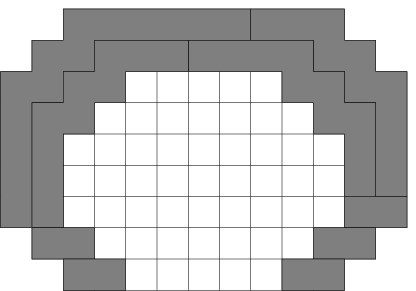}
\end{equation*}

\item After one application of the parity theorem, we can use Corollary~\ref{samelength}.
\begin{equation*}
\epsfig{file=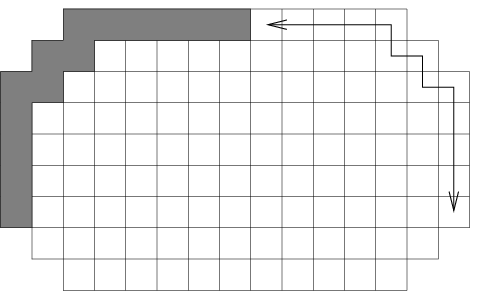}
\end{equation*}

\end{enumerate}

\end{proof}


\begin{thebibliography}{99}

\bibitem{ciucu} M.~Ciucu, Enumeration of perfect matchings in graphs with reflective symmetry, J.~Combin.~Theory Ser.~A 77 (1997) 67-97.

\bibitem{jockusch} W.~Jockusch, Perfect matchings and perfect squares, J.~Combin.~Theory Ser.~A 67 (1994) 100-115.

\bibitem{pachter} L.~Pachter, Combinatorial approaches and conjectures for $2$-divisibility problems concerning domino tilings of polyominoes, Electron.~J.~Comb. 4 (1997) R29.

\bibitem{tauraso} R.~Tauraso, A new domino tiling sequence, Journal of Integer Sequences 7 (2004) 04.2.3.

\end{thebibliography}
\end{document}